\documentclass[12pt]{article}%
\usepackage{graphicx}
\usepackage[intlimits]{amsmath}
\usepackage{latexsym}
\usepackage{amsfonts}
\usepackage{amssymb}%
 \makeatletter
\newfont{\footsc}{cmcsc10 at 8truept}
\newfont{\footbf}{cmbx10 at 8truept}
\newfont{\footrm}{cmr10 at 10truept}
\newtheorem{theorem}{Theorem}

\newtheorem{fact}[theorem]{Fact}

\newenvironment{proof}[1][Proof]{\noindent{\textbf {#1}  }}  {\hfill$\Box$\bigskip}

\begin{document}

\title{Stability for large forbidden subgraphs}
\author{Vladimir Nikiforov\\{\small Department of Mathematical Sciences, University of Memphis, Memphis,
TN 38152}\\{\small email: vnikifrv@memphis.edu}}
\maketitle

\begin{abstract}
We strengthen the stability theorem of Erd\H{o}s and Simonovits.

Write $K_{r}\left(  s_{1},\ldots,s_{r}\right)  $ for the complete $r$-partite
graph with parts of size $s_{1},\ldots,s_{r}$ and $T_{r}\left(  n\right)  $
for the $r$-partite Tur\'{a}n graph of order $n.$ Our main result is:

For all $r\geq2$ and all sufficiently small $c>0,$ $\varepsilon>0,$ every
graph $G$ of sufficiently large order $n$ with $\left\lceil \left(
1-1/r-\varepsilon\right)  n^{2}/2\right\rceil $ edges satisfies one of the conditions:

(a) $G$ contains a $K_{r+1}\left(  \left\lfloor c\ln n\right\rfloor
,\ldots,\left\lfloor c\ln n\right\rfloor ,\left\lceil n^{1-\sqrt{c}%
}\right\rceil \right)  ;$

(b) $G$ differs from $T_{r}\left(  n\right)  $ in fewer than $\left(
\varepsilon^{1/3}+c^{1/\left(  3r+3\right)  }\right)  n^{2}$ edges.\medskip

\textbf{Keywords: }\textit{stability, forbidden subgraph, }$r$\textit{-partite
subgraph; Erd\H{o}s-Simonovits theorem}

\end{abstract}

This note is part of an ongoing project aiming to renovate some classical
results in extremal graph theory, see, e.g., \cite{BoNi04} and
\cite{Nik07,Nik07c}.

Let $K_{r}\left(  s_{1},\ldots,s_{r}\right)  $ be the complete $r$-partite
graph with parts of size $s_{1},\ldots,s_{r},$ and let $T_{r}\left(  n\right)
$ be the $r$-partite Tur\'{a}n graph of order $n$. Recall the classical
stability theorem proved independently by Erd\H{o}s \cite{Erd66}, \cite{Erd68}
and Simonovits \cite{Sim68}:\medskip

\emph{For every }$r\geq2,$\emph{ }$\varepsilon>0,$\emph{ and }$\left(
r+1\right)  $\emph{-chromatic graph }$F,$\emph{ there exists }$\delta>0$\emph{
such that if a graph }$G$\emph{ of order }$n$\emph{ has }$\left\lceil \left(
1-1/r-\delta\right)  n^{2}/2\right\rceil $\emph{ edges, then either }$G$\emph{
contains }$F,$\emph{ or }$G$\emph{ differs from }$T_{r}\left(  n\right)
$\emph{ in at most }$\varepsilon n^{2}$\emph{ edges.}\medskip

Below we show that, instead of a fixed graph $F,$ we can take $\left(
r+1\right)  $-chromatic graphs whose order grows as $\log n$.

\begin{theorem}
\label{th1}Let $r\geq2,$ $1/\ln n<c<r^{-3\left(  r+14\right)  \left(
r+1\right)  },$ $0<\varepsilon<r^{-24},$ and $G$ be a graph of order $n.$ If
$\ G$ has $\left\lceil \left(  1-1/r-\delta\right)  n^{2}/2\right\rceil $
edges, then $G$ satisfies one of the conditions:

(a) $G$ contains a $K_{r+1}\left(  \left\lfloor c\ln n\right\rfloor
,\ldots,\left\lfloor c\ln n\right\rfloor ,\left\lceil n^{1-\sqrt{c}%
}\right\rceil \right)  ;$

(b) $G$ differs from $T_{r}\left(  n\right)  $ in fewer than $\left(
\varepsilon^{1/3}+c^{1/\left(  3r+3\right)  }\right)  n^{2}$ edges.
\end{theorem}

\subsubsection*{Remarks}

\begin{itemize}
\item[-] The relation between the order of the graph $F$ and the difference
between $G$ and $T_{r}\left(  n\right)  $ is hidden in the
Erd\H{o}s-Simonovits theorem. Theorem \ref{th1} makes this relation explicit.

\item[-] The relation between $c$ and $n$ in Theorem \ref{th1} needs
explanation. First, for fixed $c,$ it shows how large must be $n$ to get a
valid conclusion. But, in fact, the relation is subtler, for $c$ itself may
depend on $n,$ e.g., letting $c=1/\ln\ln n,$ the conclusion is meaningful for
sufficiently large $n.$

\item[-] Choosing randomly a graph of order $n$ with $\left\lceil \left(
1-1/r\right)  n^{2}/2\right\rceil $ edges, we can find a graph containing no
$K_{2}\left(  \left\lfloor c^{\prime}\ln n\right\rfloor ,\left\lfloor
c^{\prime}\ln n\right\rfloor \right)  $ and differing from $T_{r}\left(
n\right)  $ in more that $c^{\prime\prime}n^{2}$ edges for some positive
$c^{\prime}$ and $c^{\prime\prime},$ independent of $n$. Hence, condition
\emph{(a)} is essentially best possible.

\item[-] The factor $\varepsilon^{1/3}+c^{1/\left(  3r+3\right)  }$ in
condition \emph{(b)} is far from the best one, but is simple.
\end{itemize}

\bigskip

To prove Theorem \ref{th1}, we introduce two supporting results. Our notation
follows \cite{Bol98}; given a graph $G,$ we write:

- $\left\vert G\right\vert $ for the number of vertices set of $G;$

- $e\left(  G\right)  $ for the number of edges of $G;$

- $\delta\left(  G\right)  $ for the minimum degree of $G;$

- $k_{r}\left(  G\right)  $ for the number of $r$-cliques of $G.$\bigskip

An\emph{ }$r$\emph{-joint }of size $t$ is the union of $t$ distinct
$r$-cliques sharing an edge. We write $js_{r}\left(  G\right)  $ for the
maximum size of an $r$-joint in a graph $G.$\bigskip

The following two facts play crucial roles in our proof.

\begin{fact}
[\cite{BoNi04}, Theorem 9]\label{stabj} Let $r\geq2,$ $\ n>r^{8}$,
$\ 0<\alpha<r^{-8}/8,$ and $G$ be a graph of order $n.$ If $e\left(  G\right)
>\left(  1-1/r-\alpha\right)  n^{2},$ then either $js_{r+1}\left(  G\right)
>n^{r-1}/r^{r+6},$ or $G$ contains an induced $r$-partite subgraph $G_{0}$
satisfying $\left\vert G_{0}\right\vert \geq\left(  1-\sqrt{2\alpha}\right)
n$ and $\delta\left(  G_{0}\right)  >\left(  1-1/r-2\sqrt{2\alpha}\right)
n.\hfill\square$
\end{fact}

\begin{fact}
[\cite{Nik07}, Theorem 1]\label{ES}Let $r\geq2,$ $c^{r}\ln n\geq1,$ and $G$ be
a graph of order $n$. If $k_{r}\left(  G\right)  \geq cn^{r},$ then $G$
contains a $K_{r}\left(  s,\ldots,s,t\right)  $ with $s=\left\lfloor c^{r}\ln
n\right\rfloor $ and $t>n^{1-c^{r-1}}.\hfill\square$
\end{fact}

\bigskip

\begin{proof}
[\textbf{Proof of Theorem \ref{th1}}]Let $G$ be a graph of order $n$ with
$e\left(  G\right)  >\left(  1-1/r-\varepsilon\right)  n^{2}/2.$ Define the
procedure $\mathcal{P}$ as follows:\medskip

\textbf{While}\emph{ }$js_{r+1}\left(  G\right)  >n^{r-1}/r^{r+6}$ \textbf{do}

\qquad\emph{Select an edge contained in }$\left\lceil n^{r-1}/r^{r+6}%
\right\rceil $\emph{ cliques of order }$r+1$ \emph{and remove it from }%
$G.$\medskip

Set for short $\theta=c^{1/\left(  r+1\right)  }r^{r+6}$ and assume first that
$\mathcal{P}$ removes at least $\left\lceil \theta n^{2}\right\rceil $ edges
before stopping. Then
\[
k_{r+1}\left(  G\right)  \geq\theta n^{r-1}/r^{r+6}=c^{1/\left(  r+1\right)
}n^{r+1},
\]
and Fact \ref{ES} implies that $K_{r+1}\left(  \left\lfloor c\ln
n\right\rfloor ,\ldots,\left\lfloor c\ln n\right\rfloor ,\left\lceil
n^{1-\sqrt{c}}\right\rceil \right)  \subset G.$ Thus condition \emph{(a)}
holds, completing the proof.

Assume therefore that $\mathcal{P}$ removes fewer than $\left\lceil \theta
n^{2}\right\rceil $ edges before stopping. Writing $G^{\prime}$ for the
resulting graph, we see that
\[
e\left(  G^{\prime}\right)  >e\left(  G\right)  -\theta n^{2}>\left(
1-1/r-\varepsilon-\theta\right)  n^{2}/2.
\]

From $\ln n\geq1/c\geq r^{3\left(  r+14\right)  \left(  r+1\right)  }$ we
easily get $n>r^{8}.$ Since
\[
\varepsilon+\theta\leq r^{-24}+r^{-3\left(  r+14\right)  }r^{r+6}<r^{-8}/8
\]
and $js_{r+1}\left(  G^{\prime}\right)  <n^{r-1}/r^{r+6},$ Fact \ref{stabj}
implies that $G^{\prime}$ contains an induced $r$-partite subgraph $G_{0}$
satisfying $\left\vert G_{0}\right\vert \geq\left(  1-\sqrt{2\left(
\varepsilon+\theta\right)  }\right)  n$ and $\delta\left(  G_{0}\right)
>\left(  1-1/r-2\sqrt{2\left(  \varepsilon+\theta\right)  }\right)  n.$

Let $V_{1},\ldots,V_{r}$ be the parts of $G_{0}.$ For every $i\in\left[
r\right]  ,$ we see that%
\[
\left\vert V_{i}\right\vert \geq n-%
{\textstyle\sum\limits_{s\in\left[  r\right]  \backslash\left\{  i\right\}  }}
\left\vert V_{s}\right\vert \geq n-\left(  r-1\right)  \left(  n-\delta\left(
G_{0}\right)  \right)  \geq\left(  1/r-2\left(  r-1\right)  \sqrt{2\left(
\varepsilon+\theta\right)  }\right)  n.
\]

For each $i\in\left[  r\right]  ,$ select a set $U_{i}\subset V_{i}$ with
\[
\left\vert U_{i}\right\vert =\left\lceil \left(  1/r-2\left(  r-1\right)
\sqrt{2\left(  \varepsilon+\theta\right)  }\right)  n\right\rceil ,
\]
and write $G_{1}$ for the graph induced by $\cup_{i=1}^{r}U_{i}$. Clearly
$G_{1}$ can be made complete $r$-partite by adding at most
\[
\left(  \left(  1-1/r\right)  \left\vert G_{1}\right\vert -\delta\left(
G_{1}\right)  \right)  \left\vert G_{1}\right\vert /2
\]
edges. We see that%
\[
\delta\left(  G_{1}\right)  \geq\delta\left(  G_{0}\right)  -\left\vert
G_{0}\right\vert +\left\vert G_{1}\right\vert \geq-n/r-2\sqrt{2\left(
\varepsilon+\theta\right)  }n+\left\vert G_{1}\right\vert ,
\]
and so,
\begin{align*}
\left(  1-1/r\right)  \left\vert G_{1}\right\vert -\delta\left(  G_{1}\right)
&  \leq\left(  1/r+2\sqrt{2\left(  \varepsilon+\theta\right)  }\right)
n-\left\vert G_{1}\right\vert /r\\
&  =\left(  1/r+2\sqrt{2\left(  \varepsilon+\theta\right)  }\right)  n-\left(
1/r-2\left(  r-1\right)  \sqrt{2\left(  \varepsilon+\theta\right)  }\right)
n\\
&  =2r\sqrt{2\left(  \varepsilon+\theta\right)  }n.
\end{align*}
Therefore, $G_{1}$ can be made complete $r$-partite by adding at most
\[
2r\sqrt{2\left(  \varepsilon+\theta\right)  }n\left\vert G_{1}\right\vert
/2<r\sqrt{2\left(  \varepsilon+\theta\right)  }n^{2}%
\]
edges.

The complete $r$-partite graph with parts $U_{1},\ldots,U_{r}$ can be
transformed into $T_{r}\left(  n\right)  $ by changing at most $\left(
n-\left\vert G_{1}\right\vert \right)  n$ edges. Since
\[
\left(  n-\left\vert G_{1}\right\vert \right)  n\leq\left(  n-r\left(
1/r-\left(  r-1\right)  \sqrt{8\left(  \varepsilon+\theta\right)  }\right)
n\right)  n=2r\left(  r-1\right)  \sqrt{2\left(  \varepsilon+\theta\right)
}n^{2},
\]
we find that $G$ differs from $T_{r}\left(  n\right)  $ in at most
\[
\left(  \theta+\left(  2r^{2}-r\right)  \sqrt{2\left(  \varepsilon
+\theta\right)  }\right)  n^{2}%
\]
edges. Now, condition \emph{(b)} follows in view of
\begin{align*}
\theta+\left(  2r^{2}-r\right)  \sqrt{2\left(  \varepsilon+\theta\right)  }
&  <\theta+2\left(  2r^{2}-r\right)  \varepsilon^{1/2}+2\left(  2r^{2}%
-r\right)  \theta^{1/2}\\
&  <4r^{2}\varepsilon^{1/2}+4r^{r/2+5}c^{1/\left(  2r+2\right)  }%
<\varepsilon^{1/3}+c^{1/\left(  3r+3\right)  }.
\end{align*}
The proof is completed.
\end{proof}

\subsubsection*{Concluding remark}

Finally, here are the principles we try to follow in extremal problems:

- find results that can be used as wide-range tools, like Facts \ref{stabj}
and \ref{ES};

- give explicit conditions for the parameters in statements, like the
conditions for $\varepsilon,c,$ and $n$ in Theorem \ref{th1};

- prefer simple to optimal bounds, like the exponent $1-\sqrt{c}$ and the
factor $\varepsilon^{1/3}+c^{1/\left(  3r+3\right)  }$ in Theorem \ref{th1}.

We aim to give results that can be used further, hoping to add more integrity
to extremal graph theory.

\end{document}